\def\draft{n}
\documentclass{amsart}
\usepackage{fullpage,epic,eepic,epsfig,amssymb,pb-diagram,lamsarrow,pb-lams}

\theoremstyle{plain}

\newtheorem{theorem}{Theorem}
\newtheorem{proposition}{Proposition}[section]
\newtheorem{lemma}[proposition]{Lemma}
\newtheorem{corollary}[proposition]{Corollary}

\theoremstyle{definition}
\newtheorem{definition}[proposition]{Definition}

\theoremstyle{remark}
\newtheorem{example}[proposition]{Example}

\newtheorem{remark}[proposition]{Remark}

\def\printname#1{
	\if\draft y
		\smash{\makebox[0pt]{\hspace{-0.5in}
			\raisebox{8pt}{\tt\tiny #1}}}
	\fi
}

\newcommand{\psdraw}[2]
         {\begin{array}{c} \hspace{-1.3mm}
	\raisebox{-4pt}{\epsfig{figure=draws/#1.eps,width=#2}}
	\hspace{-1.9mm}\end{array}}

\newlength{\standardunitlength}
\catcode`\@=11
\long\def\@makecaption#1#2{%
     \vskip 10pt

\setbox\@tempboxa\hbox{
       \small\sf{\bfcaptionfont #1. }\ignorespaces #2}%
     \ifdim \wd\@tempboxa >\captionwidth {%
         \rightskip=\@captionmargin\leftskip=\@captionmargin
         \unhbox\@tempboxa\par}%
       \else
         \hbox to\hsize{\hfil\box\@tempboxa\hfil}%
     \fi}
\font\bfcaptionfont=cmssbx10 scaled \magstephalf
\newdimen\@captionmargin\@captionmargin=2\parindent
\newdimen\captionwidth\captionwidth=\hsize
\catcode`\@=12

\def\lbl#1{\label{#1}\printname{#1}}

\def\BZ{\mathbb Z}
\def\BQ{\mathbb Q}
\def\BR{\mathbb R}

\def\A{\mathcal A}

\def\cL{\mathcal L}
\def\G{\mathcal G}
\def\O{\mathcal O}

\def\F{\mathcal F}

\def\sfY{\operatorname{\mathsf{Y}}}

\def\sfI{\operatorname{\mathsf{I}}}

\def\aa{\alpha}
\def\bb{\beta}
\def\La{\Lambda}

\def\Ga{\Gamma}
\def\S{\Sigma}

\def\ga{\gamma}

\def\s{\sigma}
\def\D{\Delta}

\def\vphi{\varphi}

\def\ihs{integral homology 3-sphere}
\def\qhs{rational homology 3-sphere}

\def\fti{finite type invariant}

\def\promod{\mathop{{\rm mod}}\nolimits}
\def\AS{\mathrm{AS}}
\def\IHX{\mathrm{IHX}}
\def\ad{\mathrm{ad}}

\def\oS{\overline{S}}

\def\la{\langle}
\def\ra{\rangle}
\def\y1x{\underset{y1x}\ast}

\def\lk{\mathrm{lk}}
\def\ti{\widetilde}

\def\Fl#1{\mathcal F^{\mathrm{null}}_{#1}}

\def\GY#1{\mathcal G_{#1}(M)}
\def\Gl#1{\mathcal G^{\mathrm{null}}_{#1}}

\def\e{\varepsilon}

\def\ygraph{$\mathrm{Y}$-graph}

\def\tcircle{\circlearrowleft}
\def\tline{\uparrow}

\def\sminus{\smallsetminus}

\def\lgt{\widetilde{\mathrm{lk}}^{\ga}}

\def\fg{f^{\ga}}
\def\phig{\phi^{\ga}}

\def\longto{\longrightarrow}

\def\Lloc{\Lambda_{\mathrm{loc}}}
\def\Ql{Q^{\mathrm{loop}}}

\def\hair{\mathrm{Hair}}

\def\strutb#1#2#3{\overset{#1}{\underset{#2}{ 
\begin{array}{c} \vspace{0.0cm}
\uparrow 
\vspace{-0.25cm} \\        
| \vspace{-0.45cm} \\      
\bullet \vspace{0.00cm}   
\end{array} }}\! #3}

\def\st#1#2{\overset{#1}{\underset{#2}{
\begin{array}{c} \hspace{-1.3mm}
	\raisebox{-4pt}{\psfig{figure=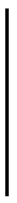,height=0.2in} }
	\hspace{-1.9mm}\end{array} }}}

\def\SO{(S^3,\O)}
\def\Zrat{Z^{\mathrm{rat}}}
\def\Aut{\mathrm{Aut}}

\def\Sym{\mathrm{Sym}}

\begin{document}


\title[The loop expansion, 
the null-move and $S$-equivalence]{
The loop expansion of the Kontsevich integral, 
the null-move and $S$-equivalence}
\author{Stavros Garoufalidis}
\address{School of Mathematics \\
          Georgia Institute of Technology \\
          Atlanta, GA 30332-0160, USA. }
\email{stavros@math.gatech.edu}
\author{Lev Rozansky}
\address{Department of mathematics \\
         University of North Carolina \\
         Chapel Hill, NC 27599}
\email{rozansky@math.unc.edu}

\thanks{The  authors are partially supported by NSF grants
        DMS-98-00703  and  DMS-97-04-893 respectively.\newline
        This and related preprints can also be obtained at
{\tt http://www.math.gatech.edu/$\sim$stavros } 
\newline
1991 {\em Mathematics Classification.} Primary 57N10. Secondary 57M25.
\newline
{\em Key words and phrases:} Kontsevich integral, claspers, null-move,
$S$-equivalence, Blanchfield pairing, Euler degree, beads, trivalent graphs,
hair map, hairy struts, hairy vortices, finite type invariants, 
$n$-equivalence.
}


\date{
June 5, 2003 \hspace{0.3cm} 
First edition: December 23, 1999.}

\begin{abstract}
The Kontsevich integral of a knot is a graph-valued invariant which (when
graded by the Vassiliev degree of graphs) is characterized by a universal 
property; namely it is a universal Vassiliev invariant of knots.
We introduce a second grading of the Kontsevich integral, the Euler degree,
and a geometric null-move on the set of knots. We explain the relation of
the null-move to $S$-equivalence, and the relation to the Euler grading
of the Kontsevich integral. The null move leads in a natural way to the
introduction of trivalent graphs with beads, and to a conjecture on a 
rational version of the Kontsevich integral, formulated by the second author
and proven in \cite{GK1}.
\end{abstract}

\maketitle



\section{Introduction}
\lbl{sec.intro}



\subsection{What is the Kontsevich integral?}
\lbl{sub.what}

The Kontsevich integral of a knot in $S^3$ is a graph-valued invariant, 
which is characterized by a universal property, namely it is a universal 
$\BQ$-valued Vassiliev invariant of knots. Using the language of physics, 
the Kontsevich integral is the Feynmann diagram expansion (i.e., perturbative 
expansion)
of the Chern-Simons path integral, expanded around a trivial flat connection,
\cite{R1,Wi}.
This explains the shape of the graphs (namely, they have univalent and 
trivalent vertices only), their vertex-orientation and their Vassiliev 
degree, namely half the number of vertices.

In physics one often takes the logarithm of a Feynmann diagram expansion,
which is given by a series of connected graphs, and one studies the
the terms with a fixed number of loops (i.e., with a fixed first betti number).

The purpose of this paper is to study the Kontsevich integral, graded
by the Euler degree (rather than the Vassiliev degree, or the number of
loops), to introduce a null-move on knots, and to relate the Euler degree
to the null-move.

This point of view explains in a natural way a Rationality Conjecture of the 
Kontsevich integral (formulated by the second author and proven by \cite{GK1}),
the relation of the null-move to $S$-equivalence, the relation of the 
null-move to cyclic branched covers \cite{GK2}, and offers an opportunity
to use Vassiliev invariants as obstructions to the existence of Seifert
forms, \cite{GT}.

\subsection{The Kontsevich integral, and its grading by the Euler
degree}
\lbl{sub.loop}

We begin by explaining the Euler expansion of the Kontsevich integral $Z$,
defined for links in $S^3$ by Kontsevich, \cite{Ko}, and extended
to an invariant of links in arbitrary closed 3-manifolds by 
the work of Le-Murakami-Ohtsuki, \cite{LMO}. Our paper will focus on
the above invariant $Z(M,K)$ of knots $K$ in \ihs s $M$:
$$
Z: \, \mathrm{Pairs}(M,K) \longto \A(\ast)
$$
where $\A(\ast)$ is the completed vector space over $\BQ$ generated 
by unitrivalent graphs with vertex-orientations modulo the well-known
antisymmetry $\AS$ and $\IHX$ relations. The graphs in question have
a Vassiliev degree (given by half the number of vertices) and the completion
of $\A(\ast)$ refers to the Vassiliev degree.

We now introduce a second grading on unitrivalent graphs: the {\em Euler
degree} $e(G)$ of a unitrivalent graph $G$ is the number of trivalent vertices
that remain when we shave-off all legs of $G$. It is easy to see that
$e(G)=-2\chi(G)$ where $\chi$ is the Euler characteristic of $G$, which 
explains the naming of this degree. The $\AS$ and $\IHX$ relations are 
homogeneous with respect to the Euler degree, thus we can let $Z_n$
denote the Euler degree $n$ part of the Kontsevich integral.

Notice that there are finitely many unitrivalent graphs with Vassiliev degree
$n$, but infinitely many connected graphs with Euler degree $n$.
For example, wheels with $n$ legs have Vassiliev degree $n$ but Euler
degree $0$.
The Kontsevich integral of a knot, graded by the Vassiliev degree is a 
universal Vassiliev \fti . On the other hand, $Z_n$ are not
Vassiliev invariants; they are rather power series of Vassiliev invariants.

Why consider the Euler degree? A deep and unexpected geometric reason is the 
content of the next section.

\subsection{The $\BZ$-null move}
\lbl{sub.nloop}

A beautiful theory of Goussarov and Habiro is the study of the geometric notion
of {\em surgery on a clasper}, see \cite{Gu1, Gu2, Ha} and \cite{GGP} 
(in the latter
claspers were called clovers). Following the notation of \cite{GGP}, given
a clasper $G$ in a 3-manifold $N$, we let $N_G$ denote the result of surgery.
Clasper surgery can be described in terms of twisting genus $3$ handlebodies
in $N$, or alternatively in terms of surgery on a framed link in $N$. We
will refer the reader to \cite{GGP} for the definition and conventions of
surgery on a clasper. We will {\em exclude} claspers with no trivalent
vertices, that is, ones that generate $I$-moves.

In the present paper we are interested in pairs $(M,K)$ and claspers $G
\subset M\sminus K$ whose leaves are {\em null homologous links} in 
$M\sminus K$. We will call such claspers $\BZ$-{\em null}, or simply
{\em null}. The terminology is motivated by the fact that the leaves of such 
claspers are sent to $0$ under the map: $\pi_1(M\sminus K)\to H_1(M\sminus K)
\cong \BZ$.

In order to motivate our interest in null claspers, recall the basic and
fundamental principle: surgery on a clasper preserves the homology and (when
defined) the linking form. Applying that principle in the $\BZ$-cover of
a knot complement, it follows that surgery on a null clasper preserves
the Blanchfield form, as we will see below.

Surgery on null claspers describes a move on the set of knots in \ihs s.
Given this move, one can define in the usual fashion a notion of {\em \fti s}
and a dual notion of $n$-{\em equivalence}, as was explained by Goussarov and 
Habiro. Explicitly, we can consider a decreasing filtration $\Fl {}$ on the
vector space generated by all pairs $(M,K)$ as follows: $\Fl n$ is generated
by $[(M,K),G]$, where $G=\{G_1,\dots, G_n\}$ is a disjoint collection of null 
claspers in $M\sminus K$, and where
\begin{equation}
\lbl{eq.MKG}
[(M,K),G]=\sum_{I \subset \{0,1\}^n} (-1)^{|I|} (M,K)_{G_I}
\end{equation}
(where $|I|$ denotes the number of elements of $I$ and $N_{G_I}$ stands 
for the result of simultaneous surgery on claspers $G_i \subset N$ for all 
$i \in I$).

\begin{definition}
\lbl{def.fti}
A function $f: \, \mathrm{Pairs}(M,K) \longto \BQ$
is a {\em \fti\ of null-type} $n$ iff $f(\Fl {n+1})=0$.
\end{definition}

We will discuss the general notion of $n$-null-equivalence in a later 
publication, at present we will consider the special case of $n=0$. 

\begin{definition}
\lbl{def.0equiv}
Two pairs $(M,K)$ and $(M',K')$ are {\em null-equivalent} iff one can
be obtained from the other by a sequence of null moves.
\end{definition}
 
In order to motivate the next lemma, recall a result of 
Matveev \cite{Ma}, who showed that two closed 3-manifolds $M$ and $M'$ are 
equivalent under a sequence of clasper surgeries iff they have the same 
homology and linking form. This answer is particularly pleasing since it is 
expressed in terms of abelian algebraic invariants.

\begin{lemma}
\lbl{lem.deg0}
The following are equivalent: $(M,K)$ and $(M',K')$
\begin{itemize}
\item[(a)]
are null equivalent 
\item[(b)]
are $S$-equivalent 
\item[(c)]
have isometric Blanchfield form.
\end{itemize}
\end{lemma}

\begin{proof}
For a definition of $S$-{\em equivalence} and {\em Blanchfield forms}, see 
for example
\cite{Ka,L2}. It is well-known that $(M,K)$ and $(M',K')$ are $S$-equivalent
iff they have isometric Blanchfield forms; for an algebraic-topological
proof of that combine Levine and Kearton, \cite{L1,Ke}, or alternatively
Trotter, \cite{Tr}. Thus (b) is equivalent to (c).

If $(M',K')$ is obtained from $(M,K)$ by a null move on a clasper $G$,
then $G$ lifts to the universal abelian cover $\ti X$ of the knot complement
$X=M\sminus K$. The lift $\ti G$ of $G$ is a union of claspers translated by
by the group $\BZ$ of deck transformations. Furthermore, $\ti {X_G}$ is 
obtained from $\ti X$ by surgery
on $\ti G$. Since clasper surgery preserves the homology and linking form, it
follows that $(M',K')$ and $(M,K)$ have isometric Blanchfield forms,
thus (a) implies (c).

Suppose now that $(M',K')$ is $S$-equivalent to $(M,K)$. Then, by Matveev's
result, the  \ihs\ $M$ can be obtained from
$S^3$ by surgery on some null $G$. Moreover, if $K$ is a knot in $M$ 
that bounds a Seifert surface $\Sigma$, we can always arrange (by an
isotopy on $G$) that $G$ is disjoint from $\Sigma$, thus $G$ is $(M,K)$-null.
Thus, if $(M,K)$ and $(M',K')$ are $S$-equivalent pairs,
modulo null moves we can assume that $M=M'=S^3$. In that case,
Naik and Stanford show that $S$-equivalent pairs $(S^3,K)$ and $(S^3,K')$
are equivalent under a sequence of double $\D$elta-moves, that is 
null-moves where all the leaves of the clasper bound disjoint disks, 
\cite{NS}. Thus, (b) implies (a) and the Lemma follows.
\end{proof}

It follows that if $(M,K)$ and $(M',K')$ are null 
equivalent, then they have the same Alexander module, in particular the same
Alexander polynomial, and the same
algebraic concordance invariants which were classified by Levine 
\cite{L1,L2}.

\begin{remark}
\lbl{rem.linkingform}
There is a well-known similarity between $\BQ/\BZ$-valued linking forms
of rational homology spheres and $\BQ(t)/\BZ[t,t^{-1}]$-valued Blanchfield
forms. From this point of view,  \ihs\ correspond to
knots with trivial Alexander polynomial. This classical analogy extends to 
the world of \fti s as follows:
one can consider \fti s of \qhs s (using surgery on claspers) and 
\fti s of knots in \ihs s (using surgery on null-claspers). In our paper,
we will translate well-known results about \fti s and $n$-equivalence of
\ihs s to the setting of knots with trivial Alexander polynomial.
\end{remark}

Are there any nontrivial \fti s of null-type? The next result answers 
this, and reveals an unexpected relation between the null-move and the
Euler degree:

\begin{theorem}
\lbl{thm.loop2}
For all $n$, $Z_{n}$ is
a null-type $n$ invariant of pairs $(M,K)$ with values in $\A_{n}(\ast)$. 
\end{theorem}

In the above, $\A_n(\ast)$ stands for the subspace of $\A(\ast)$ 
generated by diagrams of Euler-degree $n$.

The above theorem has implications for the {\em loop expansion} of the
Kontsevich integral defined as follows. The logarithm of the Kontsevich
integral 
$$
\log Z: \, \mathrm{Pairs}(M,K) \longto \A^c(\ast)
$$
takes values in the completed  vector space (with respect to the Vassiliev 
degree) $\A^c(\ast)$ generated by connected vertex-oriented unitrivalent 
graphs, modulo the $\AS$ and $\IHX$ relations. Let us define the {\em loop
degree} of a graph to be the number of loops, i.e., the first betti number.
The $\AS$ and $\IHX$ relations are homogeneous with respect to the loop
degree, thus we can define $\Ql_n$ to be the $(n+1)$-loop degree of
$\log Z$. Since a connected trivalent graph with $n+1$ loops has Euler degree
$2n$, and since $\log Z$ is an additive invariant of pairs under connected
sum, it follows that 

\begin{corollary}
\lbl{cor.loop2}
For all $n$, $\Ql_n$ is an additive null-type $2n$ invariant of pairs $(M,K)$
with values in the loop-degree $(n+1)$ part of $\A^c(\ast)$.
\end{corollary}

Our next task is to study the graded quotients $\Gl n=
\Fl n/\Fl {n+1}$; knowing this would essentially answer the question of
how many null-\fti s there are.

Since null equivalence has more than one equivalence classes, we ought to
study $\Fl {}$ and $\Gl {}$ on each null equivalence class. This is formalized
in the following way. Given a pair $(M,K)$, let $\Fl {}(M,K)$ denote the
subspace of $\Fl {}$ generated by all pairs $(M',K')$ which are null
equivalent to $(M,K)$. It is easy to see that $\Fl {}(M,K)=\Fl {}(M',K')$
when $(M,K)$ and $(M',K')$ are null equivalent and that we have
a direct sum decomposition
$$
\Fl {}=\oplus_{0-\text{null equiv. classes} (M,K)} \Fl {}(M,K).
$$
Fixing $(M,K)$, we can define $\Gl n (M,K)=\Fl n (M,K)/\Fl {n+1} (M,K)$ 
accordingly.
The simplest case to study are the graded quotients $\Gl {}(S^3,\O)$ for the
unknot $\O$ in $S^3$. This is equivalent to the study of null-\fti\ of
knots with trivial Alexander polynomial. 

It turns out in a perhaps unexpected way that $\Gl {}(S^3,\O)$ can be
described in terms of trivalent graphs with {\em beads}. Let us define
these here. Let $\La=\BZ[t,t^{-1}]$ denote the group-ring of the integers.
It is a ring with involution $r \to \bar r$ given by $\bar t=t^{-1}$
and augmentation map $\e: \La \to \BZ$ given by $\e(t)=1$.
Consider a trivalent graph $G$ with oriented edges. A {\em $\La$-coloring} of
$G$ is an assignment of an element of $\La$ to every edge of $G$. We will 
call the assignment of an edge $e$, the {\em bead} of $e$.
 
\begin{definition}
\lbl{def.Ala}
$\A(\La)$ is the completed vector space over $\BQ$ generated by pairs
$(G,c)$ (where $G$ is a trivalent graph, with oriented edges and 
vertex-orientation and $c$ is a $\La$-coloring of $G$), modulo the
relations: $\AS$, $\IHX$, Linearity, Orientation Reversal,
Holonomy and Graph Automorphisms. 
\end{definition}

\begin{figure}[htpb]
$$ 
\psdraw{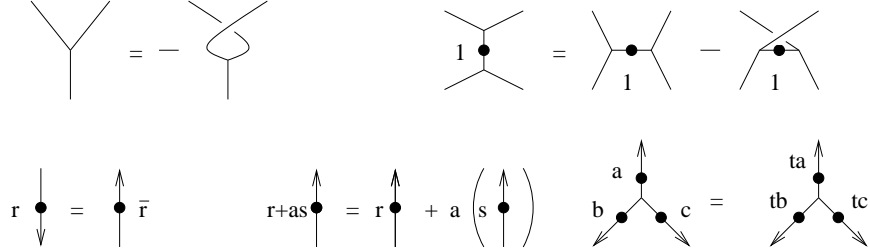}{4.5in} 
$$
\caption{The $\AS$, $\IHX$ (for arbitrary orientations of the edges),
Orientation Reversal, Linearity and Holonomy Relations.}\lbl{relations4}
\end{figure}

$\A(\La)$ is graded by the Euler degree of a graph, and the completion
is with respect to the above grading. 

\begin{remark}
\lbl{rem.integers}
We could have defined $\Gl {}(M,K)$ and $\A(\La)$ over $\BZ$ rather than over 
$\BQ$. We will not use special notation to indicate this; instead when
we wish, we will state explicitly which coefficients we are using.
\end{remark}

\begin{theorem}
\lbl{thm.loop1}
Over $\BZ$,  there is degree-preserving map:
$$
\A(\La)\to \Gl {}(S^3,\O)
$$
which is onto, over $\BZ[1/2]$. Since $\A_{\mathrm{odd}}(\La)=0$, it follows 
that over $\BZ[1/2]$,  $\Fl n(S^3,\O)$ is a $2$-step filtration, i.e., 
satisfies $\Fl {2n+1}(S^3,\O)=\Fl {2n+2}(S^3,\O)$ for all $n$. 
\end{theorem}

\subsection{Further developments}
\lbl{sub.further}

In view of Theorems \ref{thm.loop2} and \ref{thm.loop1}, it is natural to
ask if there is any relation among the algebras $\A(\ast)$ and $\A(\La)$.
In order to make contact with later work, we mention here an important
{\em hair} map
$$
\hair: \A(\La) \longto \A(\ast)
$$
which is defined by replacing a bead $t$ by an exponential of hair:
$$
\strutb{}{}{t} \to \sum_{n=0}^\infty \frac{1}{n!} 
\, \, \psdraw{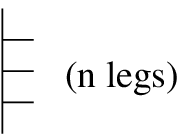}{.7in}
$$

The Kontsevich integral (graded by the Euler degree, and evaluated on
Alexander polynomial $1$ knots) fails to be a universal 
\fti\ (with respect to the null-move) because it takes values in the 
``wrong'' space, namely $\A(\ast)$ rather than $\A(\La)$. It is natural to
conjecture the existence of an invariant $\Zrat$ that fits in a 
commutative diagram
$$
\divide\dgARROWLENGTH by1
\begin{diagram}
\node[2]{\A(\La)}       
\arrow{s,r}{\hair}                     \\
\node{(\Delta=1)-\mathrm{knots}}
\arrow{e,t}{Z}\arrow{ne,t}{\Zrat}
\node{\A(\ast)}
\end{diagram}
$$
This Rationality Conjecture was formulated by the second author (in a version
for all knots) and was proven in \cite{GK1}. 

In an earlier version of the paper, we formulated a conjectural relation
of the 2-loop part of the $\Zrat$ invariant and the Casson-Walker invariant
of cyclic branched coverings of a knot. This relation has been settled by
\cite[Corollary 1.4]{GK2}.

\begin{remark}
In the study of knot theory via surgery,
knots with Alexander trivial polynomial are topologically (but not smoothly)
slice, \cite{FQ} and also \cite[Appendix]{GT}. The ``abelian'' and ''solvable''
 invariants 
of knots with trivial Alexander polynomial, such as their Alexander module, 
Casson-Gordon invariants and all of their recent systematic obstructions
of Cochran-Orr-Teichner \cite{COT} vanish. However, already the Euler degree
$2$ part of the Kontsevich integral $Z_2$ (which incidentally equals to
the $2$-loop part of $Z$) does not vanish on Alexander 
polynomial $1$ knots. The symbol of $Z_2$ can be computed in terms of
equivariant linking numbers (see Theorem \ref{thm.u2})
and this gives good realization properties for
$Z_2$. The $Z_2$ invariant offers an opportunity to settle an error in one of 
M. Freedman's lemmas on knots with Alexander polynomial $1$, \cite{GT}. 
\end{remark}

\begin{remark}
\lbl{rem.loopfti}
To those who prefer moves that untie knots, it might sound disappointing 
that the null-move fails to do so. On the other hand, the geometric null-move
describes a natural equivalence relation on knots, namely isometry of
Blanchfield forms. Concordance is another well-prized equivalence relation
on knots. If one could generate concordance in terms of surgery of some
type of claspers, this could open the door for constructing a plethora
of concordance invariants of knots. It is known that surgery on certain
claspers preserve concordance, \cite{GL2,CT}, but it is also known
that these moves do not generate concordance, \cite{GL2}.
\end{remark}

\subsection{Plan of the paper}
\lbl{sub.plan}

The paper consists of five, largely independent, sections. 

\noindent
$\bullet$
In Section \ref{sec.loop} we apply the topological calculus
of claspers to the case of null-claspers which leads naturally to trivalent
graphs with beads, and their relation to the   
graded quotients $\G (S^3,\O)$.

\noindent
$\bullet$
In Section \ref{sub.Zbrief} we give a detailed study of the behavior of
the Aarhus integral under surgery on claspers. In particular, 
counting arguments above the critical degree imply that the Euler degree $n$ 
part of the Kontsevich integral is a \fti\ of null-type $n$. 
This provides a conceptual relation between the 
null-move and the Euler degree of graphs. This counting leads naturally
to the study of the hairy strut part of the Kontsevich integral. Struts
correspond to linking numbers and the goal in the rest of the sections
is to show that hairy struts correspond to equivariant linking numbers.

\noindent
$\bullet$
In Section \ref{sec.abelian} we review linking numbers and equivariant
linking numbers, and give an axiomatic description of the latter, motivated
by the calculus of claspers.

\noindent
$\bullet$
In the final Section \ref{sub.Kbehavior} we show that the 
hairy strut part of the Kontsevich integral satisfies the same axioms
as the equivariant linking numbers, and as a consequence of a uniqueness
result, the two are equal.

\tableofcontents

\section{Surgery on claspers and the the null-filtration}
\lbl{sec.loop}

\subsection{A brief review of surgery on claspers}
\lbl{sub.cbrief}

In this section we recall briefly the definition of surgery on claspers,
with the notation of \cite{GGP}. We refer the reader to \cite{GGP} for
a detailed description.
A framed graph $G$ in a 3-manifold $M$ is called a {\em clasper of degree} $1$
(or simply, a {\em Y-graph\/}), if it is the image of $\Ga$ under a smooth 
embedding $\phi_G:N\to M$ of a neighborhood $N$ of $\Ga$. The embedding 
$\phi_G$ can be recovered from $G$ up to isotopy. Surgery on $G$ can be 
described either by surgery on the corresponding framed six component link, 
or in terms of cutting a neighborhood of $G$, (which is a handlebody of genus 
$3$), twisting by a fixed diffeomorphism of its boundary and gluing back.

\begin{figure}[htpb]
$$
\psdraw{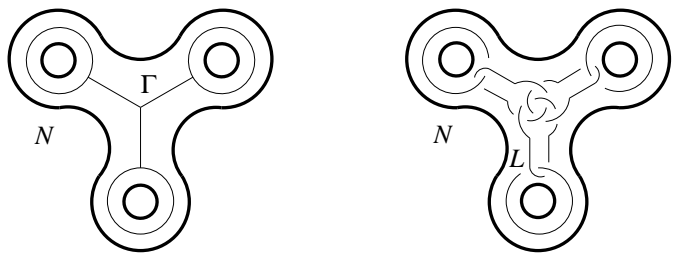}{3.5in}
$$
\caption{Y-graph and the corresponding surgery link}\lbl{yvaria}
\end{figure}

A clasper of higher degree is a (thickening of) an embedded trivalent graph,
with some distinguished loops which we call {\em leaves}. The degree
of a clasper is the number of trivalent vertices, excluding those of the 
leaves. Surgery on claspers of higher degree is defined similarly. 
For example,
a clasper of degree $2$ and its corresponding surgery link are shown below:
$$
\psdraw{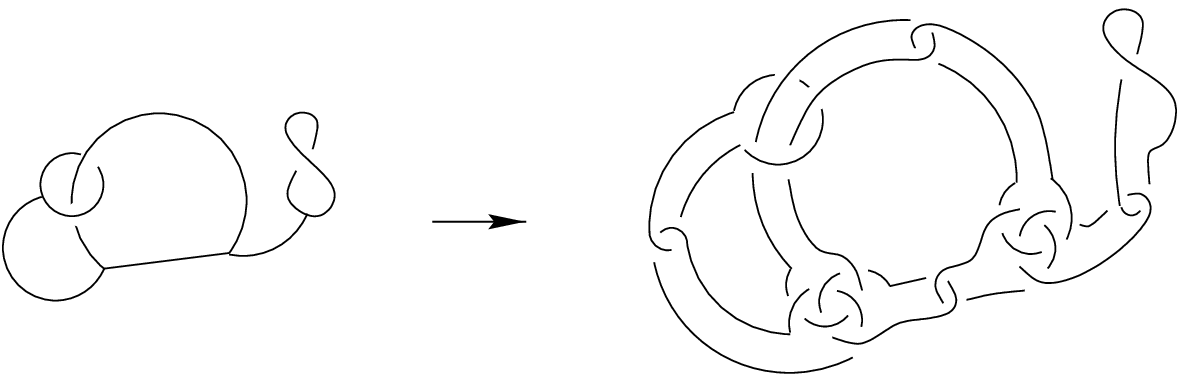}{3.5in}
$$

Surgery on a clasper of degree $n$ can be described in terms of surgery
on $n$ claspers of degree $1$. We will exclude claspers of degree $0$, that
is with no trivalent vertices. These correspond to $I$-moves in the language
of Goussarov and Habiro.

\subsection{A review of \fti s of \ihs s}
\lbl{sub.reviewGH}

As we mentioned in the introduction (see Remark \ref{rem.linkingform}), the 
study of the graded 
quotients $\G (S^3,\O)$ is entirely analogous to the 
study of the graded quotients $\GY {}(S^3)$ of \fti s of \ihs s. For a 
detailed description of the latter invariants, we refer the reader to 
\cite{Ga2}. Let us recall here some key ideas, for later use.

Let $\F(S^3)$ denote the vector space over $\BQ$ generated by \ihs s,
and let $\F_n(S^3)$ denote the subspace of $\F(S^3)$ generated by
$[M,G]$ for all disjoint collections of claspers $G=\{G_1,\dots,G_n\}$ in 
\ihs s $M$, where
$$
[M,G]=\sum_{I \subset \{0,1\}^n} (-1)^{|I|} M_{G_I}
$$
Compare with Equation \eqref{eq.MKG}.

\begin{definition}
\lbl{def.ihsfti}
A function $f: \, \text{\ihs s} \longto \BQ$
is a {\em \fti\ of type} $n$ iff $f(\F_{n+1}(S^3))=0$.
\end{definition}

Let $\G_n(S^3)=\F_n(S^3)/\F_{n+1}
(S^3)$ denote the graded quotients. These quotients can be described
in terms of a vector space $\A(\phi)$ defined as follows:

\begin{definition}
\lbl{def.Aphi}
$\A(\phi)$ is the completed vector space over $\BQ$ generated by trivalent 
graphs with vertex-orientations, modulo the $\AS$ and $\IHX$ relations.
\end{definition}

In \cite{GGP} we defined a degree-preserving map
\begin{equation}
\lbl{eq.AS3}
\A(\phi) \longto \G(S^3)
\end{equation}
as follows: Given an abstract vertex-oriented trivalent graph $G$, choose
an embedding $\phi:G \to S^3$ that preserves its vertex-orientation.
Recall that a vertex-orientation is a cyclic order of the set of three flags
around a trivalent vertex. An embedding $\phi$ gives rise at each vertex
$\phi(v)$ to a frame that consists of the tangent vectors of the images under
$\phi$ of the three flags around $v$. Comparing this frame with the standard
orientation of $S^3$, we get a number $\e_{\phi,v} = \pm 1$ at each vertex
$v$ of $G$. We say that $\phi$ is orientation preserving if $\prod_v
\e_{\phi,v}=1$.

We will consider the clasper $\phi (G) \subset S^3$ of degree $n$
and the associated element
$[S^3, \phi (G)] \in \G_n (S^3)$. Lemma \ref{lem.slide} implies that the 
element $[S^3, \phi (G)] \in \G (S^3)$ is independent of the embedding $\phi$.
Further, it was shown in \cite[Corollary 4.6]{GGP} and 
\cite[Theorem 4.11]{GGP} that the $\AS$ and the $\IHX$ relations hold.
This defines the map \eqref{eq.AS3} over $\BZ$, which is obviously 
degree-preserving.

We now recall the following theorem of \cite[Theorem 4.13]{GGP}:

\begin{theorem}
\lbl{thm.GGPAS3}
The map \eqref{eq.AS3} is onto, over $\BZ[1/2]$.
\end{theorem}

Let us comment on the theorem. It is not obvious that the map \eqref{eq.AS3} 
is onto, 
since $\G(S^3)$ is generated by  elements of the form $[M,G]$ for claspers
$G$ of degree $n$ with leaves {\em arbitrarily} framed links, whereas
the image of the map \eqref{eq.AS3} involves claspers whose leaves are linked 
in the pattern of $0$-framed Hopf links.

Let us briefly review a proof of Theorem \ref{thm.GGPAS3}, taken from
\cite{Ga2}, which adapts well to our later needs. Observe that each leaf of a 
clasper in $S^3$ is {\em null homologous}. 
Using Lemmas \ref{lem.cut} and Lemma 
\ref{lem.slide} we can write $[S^3,G]$ as a linear combination of $[S^3,G']$
such that each $G'$ satisfies the following properties:
\begin{itemize}
\item
Each leaf $l$ of $G'$ is an unknot with framing $0$ or $\pm 1$ and
bounds a disk $D_l$.
\item
The disks $D_l$ of $\pm 1$-framed leaves $l$ are disjoint from each other,
disjoint from $G'$, and disjoint from the disks of the $0$-framed leaves.
\item
Each disk of a $0$-framed leaf intersects precisely one other disk
of a $0$-framed leaf in a single clasp.
\end{itemize}
Using \cite[Lemma 4.8]{GGP}, and working over $\BZ[1/2]$, we may assume
that $G'$ as above has no $\pm 1$-framed leaves. Then, it follows that 
$[S^3,G']$ lies in the image of \eqref{eq.AS3}. This completes the proof
of Theorem \ref{thm.GGPAS3}.
\qed

\begin{remark}
\lbl{rem.qhs}
If we fix a \qhs\ $M$, we could define a filtration on the vector space
$\F(M)$ generated by all homology spheres with a fixed linking form
(these are exactly $0$-equivalent to $M$ under a sequence of clasper moves) 
and introduce a decreasing filtration on $\F(M)$ and graded quotients $\G(M)$. There is a degree-preserving map
$$
\A(\phi)\longto \G(M)
$$
which is onto, over $\BZ[1/2,1/|H_1(M)|]$. The proof uses the same reasoning
as above together with the fact that every knot (such as a leaf of a clasper)
in $M$ is null homologous in $M$, once multipled by $|H_1(M)|$.
\end{remark}

\subsection{The ``Cutting'' and ``Sliding'' Lemmas}
\lbl{sub.twolemmas}

In this section we recall the key Cutting and Sliding Lemmas.

\begin{lemma}\cite{Gu2,Ha}{\bf (Cutting a Leaf)}
\lbl{lem.cut}
Let $\Ga,\Ga'$ and $\Ga''$ be the claspers of Figure \ref{split3} in a 
handlebody $V$ embedded in a \ihs\ $M$ and let $\Ga_0$ be a clasper of 
degree $m-1$ in the complement of $V$.  If $G=\Ga_0 \cup \Ga$,
$G'=\Ga_0 \cup \Ga'$ and $G''=\Ga_0 \cup \Ga''$ 
then we have that
$$[M,G]=[M,G']+[M,G''] \,\, 
\promod \F_{m+1}(S^3).$$
\end{lemma}

\begin{figure}[htpb]
$$ 
\psdraw{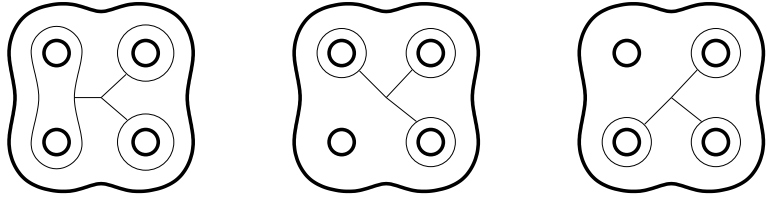}{4in} 
$$
\caption{Cutting a leaf.
The claspers $\Ga, \Ga'$ and $\Ga''$, from left to right.}\lbl{split3}
\end{figure}

Informally, we may think that we are splitting a leaf of $G$ into a connected
sum.

\begin{lemma}\cite{Gu2,Ha}{\bf (Sliding an Edge)}
\lbl{lem.slide} 
Let $\Ga^s$ and $\Ga$ be the claspers of Figure \ref{slide} in a 
handlebody $V$ embedded in a manifold $M$ and let $\Ga_0$ be a clasper of 
degree $m-1$ in the complement of $V$.  If 
$G^s=\Ga_0 \cup \Ga^s$ and $G=\Ga_0 \cup \Ga$, then 
we have that
$$[M,G^s]=[M,G] \,\, 
\promod \F_{m+1}(S^3).$$
\end{lemma}

\begin{figure}[htpb]
$$ 
\psdraw{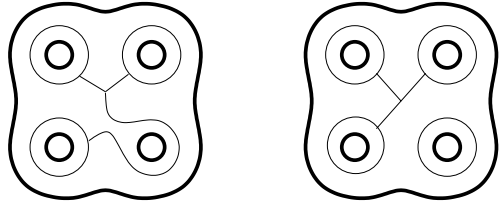}{2.6in} 
$$
\caption{Sliding an edge. The claspers $\Ga^s$ and $\Ga$.}\lbl{slide}
\end{figure}

Informally, we may think that we are sliding an edge of $G$.

\subsection{The graded quotients $\Gl {}(S^3,\O)$}
\lbl{sub.graded}

Our first goal in this section is to define the analogue of the map
\eqref{eq.AS3} using trivalent graphs with beads.
This will be achieved by introducing {\em finger moves} of $\SO$-null
claspers.

Let us begin with an alternative description of the space $\A(\La)$.
Consider a unitrivalent graph $G$ with vertex-orientation and edge 
orientation, and cut each edge to
a pair of {\em flags} (i.e., half-edges). Orient the two flags $\{e_b,e_t\}$
incident to an edge $e$ of $G$ as follows
$$
\psdraw{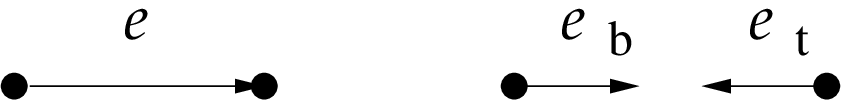}{2in}
$$

\begin{definition}
\lbl{def.LG}
The ring $\La_G$ is the polynomial ring over $\BZ$ with generators
$t_e^{\pm 1}$ for every flag $e$ of $G$ and relations:
$t_{e_b} t_{e_t}=1$ if the pair of flags $(e_b,e_t)$ is an edge of $G$,
and $\prod_{v \in e} t_e=1$ for all trivalent vertices $v$ of $G$.
\end{definition}

Let us define an abelian group:
\begin{equation}
\lbl{eq.AR}
\A^C(\La)=\left(\oplus_{G} \La_G \cdot G \right)/(\AS,\IHX,\Aut)
\end{equation}
where the sum is over all isomorphism classes of trivalent graphs with
oriented edges and vertex-orientation.

\begin{example}
\lbl{ex.Lrings}
For a strut $\sfI$, a vortex $\sfY$ and a Theta graph $\Theta$, the associated
rings are given by
\begin{eqnarray*}
\La_{\sfI} & \cong & \BZ[t^{\pm 1}] \\
\La_{\sfY} & \cong & \BZ[t_1^{\pm 1},t_2^{\pm 1},t_3^{\pm 1}]/(t_1t_2t_3-1) \\
\La_{\Theta} & \cong & \BZ[t_1^{\pm 1},t_2^{\pm 1},t_3^{\pm 1}]/(t_1t_2t_3-1,
\Sym_2 \times \Sym_3)
\end{eqnarray*} 
where $\Sym_3$ acts as a permutation on the $t_i$s and $\Sym_2$ acts
as an simultaneous involution of the powers of the $t_i$. 
\end{example}

\begin{lemma}
\lbl{lem.finger} 
Consider an $\SO$-null clasper $G$ of degree $n$, and let $G^{nl}$ 
(the superscipt stands for 'no leaves') denote
the abstract unitrivalent graph obtained from removing the leaves of
$G$. Then, there exists a geometric action
$$
\La_{G^{nl}} \times [\SO, G] \longto \G \SO
$$
given by {\em finger moves}.  
\end{lemma}

\begin{proof}
The action will be defined in terms of $I$-moves.
Let us recall surgery on a clasper with no trivalent vertices, a so-called
{\em $I$-move}:
$$
\psdraw{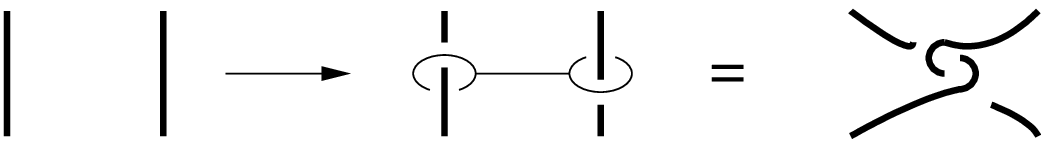}{3in}
$$
An $I$-move can be thought of as a right-handed {\em finger move}, or a 
right-handed Dehn twist, and has an inverse given by a left-handed finger
move (indicated by a stroke in the opposite direction in \cite{Gu2}).
We will call right-handed $I$-moves positive, and left-handed ones negative.

To define the action, without loss of generality, we will assume that $G$ is 
a collection of claspers of degree $1$. Further, we orient the flags of $G$
outward, towards their univalent vertices.

Consider a monomial $c=\prod_e t_e^{n(e)}$ where the product is over
the set of flags of $G$. For each flag $e$ of $G$ choose a segment $a_e$ 
of $\O$, such that segments of different flags are nonintersecting.
For each flag $e$ of $G$, choose a collection $I_{e,a_e}$ of $n(e)$ 
$I$-claspers such that when $n(e)=1$, we have:
$$
\psdraw{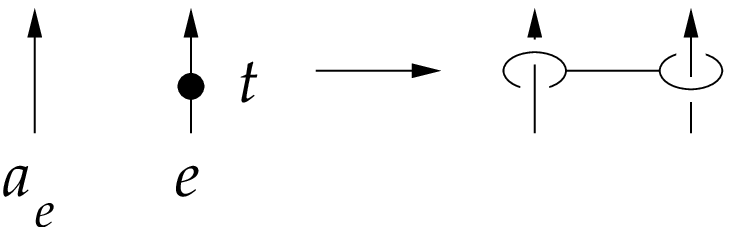}{2in} \hspace{1cm} \psdraw{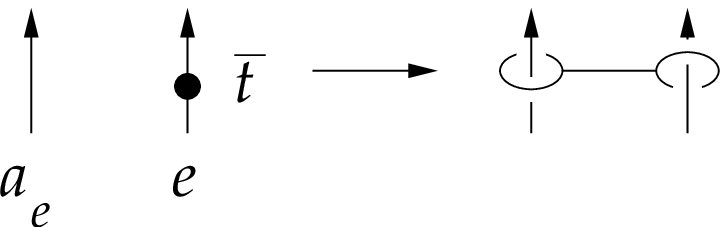}{2in}
$$
Let $I_c=\cup_e I_{e, a_e}$ be the collection of $I$-claspers for all flags
$e$ of $G$, and let $I_c \cdot G$ denote the image of $G$
after clasper surgery on $I_c$. This defines an element $[\SO, I_c\cdot G] 
\in \G \SO$.
We claim that this element depends on $c$ alone, and not on the intermediate
choices of arcs and the $I$-claspers. Indeed, Lemma \ref{lem.slide} implies
that if an edge of $I_c. G$ passes through two oppositely oriented arcs of 
$\O$, resulting in a clasper $G'$, then $[\SO,I_c\cdot G]=[\SO, G'] \in 
\G \SO$. 
This implies easily our previous claim.

Thus, we can define the action of $c$ on $G$ by $c\cdot [\SO,G]=
[\SO,I_c \cdot G]$.
It is easy to see that the relations of the ring $\La_{G^{nl}}$ hold.
\end{proof}

We can now define map 
\begin{equation}
\lbl{eq.ARla}
\A^C(\La)\to \Gl {}\SO
\end{equation}
over $\BZ$ as follows (the superscipt stands for 'clasper'): consider a 
trivalent graph $G$ with a vertex-orientation
and an orientation of the edges, and let $\phi: G \to B^3 \subset \SO$ be an
orientation preserving embedding in a small ball $B^3$. Consider the element
$[S^3,\phi(G)] \in \G \SO$, as in the map \eqref{eq.AS3}.
Using Lemma \ref{lem.finger} and the fact that $G$ is trivalent, we get
a map
$$
\La_G \longto \G \SO
$$
induced by the action $ \La_G \times [S^3,\phi(G)] \to \G \SO$. 
Just as in the map \eqref{eq.AS3}, the $\AS$, $\IHX$ and $\Aut$ relations
are preserved. This defines the map \eqref{eq.ARla}.

We now give a version of the ring $\La_G$ that uses half the number
of generators:

\begin{definition}
\lbl{def.Lahalf}
If $G$ is a vertex-oriented graph with oriented edges, let us define the 
ring $\La^V_G$ over
$\BZ$ with generators $t_e^{\pm 1}$ for every edge $e$ of $G$ and relations:
$\prod_{v \in e} t_e=1$ for all trivalent vertices $v$ of $G$.
\end{definition}

The next lemma identifies $\La_G$ with $\La^V_G$:

\begin{lemma}
\lbl{lem.Lahalf}
If $G$ has oriented edges, then we have a canonical isomorphism:
$$
\La_G \cong \La^V_G.
$$
\end{lemma}

\begin{proof}
If $e$ is an oriented edge, let $(e_b,e_t)$ denote the pair of flags such
that $e$ has the same orientation as $e_b$ and opposite from $e_t$.
This gives rise to a map $t_{e_b} \in \La_G \to t_e \in \La^V_G$ and
$t_{e_t} \in \La_G \to t_e^{-1} \in \La^V_G$, which is an isomorphism.
\end{proof}

\begin{lemma}
\lbl{lem.Aid1}
We have a canonical isomorphism 
$$
\A(\La) \cong \A^C(\La)
$$
\end{lemma}

\begin{proof}
It follows immediately using Lemma \ref{lem.Lahalf} above.
\end{proof}
 
Combining the above lemma with the map \eqref{eq.ARla}, we get the
desired map
$$
\A(\La)\to \Gl {}(S^3,\O)
$$
of Theorem \ref{thm.loop1}.

\begin{proof}(of Theorem \ref{thm.loop1})
We now claim that the proof of Theorem \ref{thm.GGPAS3} works without
change, and proves Theorem \ref{thm.loop1}. Indeed, Lemmas \ref{lem.cut}
and \ref{lem.slide} work for $\SO$-null claspers as stated.
Furthermore, if $G \subset S^3\sminus\O$ is an $\SO$-null clasper, then
each leaf of $G$ lies in the commutator group $[\pi,\pi]=1$, where
$\pi=\pi_1(S^3\sminus\O) \cong \BZ$. Since 
$1=[\pi,\pi]=[[\pi,\pi],[\pi,\pi]]$, this implies
that each leaf of $G$ bounds a surface in $S^3\sminus\O$ whose bands are 
null homologous links in $S^3\sminus\O$. The surfaces for different
leaves may intersect each other. We can apply now the proof 
of Theorem \ref{thm.GGPAS3} exactly as stated, to conclude the proof of 
Theorem \ref{thm.loop1}.
\end{proof}

We end this section with an alternative description of the ring $\La_G$,
and thus of the vector space of graphs $\A(\La)$, that was introduced
by the second author \cite{R3} and studied by P. Vogel in unpublished work.
We thank P. Vogel for explaining us his unpublished work.

\begin{lemma}
\lbl{lem.RV}
For a trivalent graph $G$ with oriented edges we have 
a canonical isomorphism
$$
\La^V_G \cong \BZ[H^1(G,\BZ)]
$$
\end{lemma}

\begin{proof}
Using Lemma \ref{lem.Lahalf}, we will rather describe a canonical
isomorphism $\La^V_G \cong \BZ[H^1(G,\BZ)]$. Recall the exact sequence
$$
0 \longto C^0(G,\BZ) \longto C^1(G, \BZ) \longto H^1(G,\BZ) \longto 0
$$
where $C^0(G,\BZ)$ and $C^1(G,\BZ)$ are the abelian groups of $\BZ$-valued
functions on the vertices and oriented edges of $G$. Let $c \in C^1(G,\BZ)$.
This gives rise to the element $\prod_e t_{e_b}^{c(e)} \in \La_G$
where the product is taken over all edges of $G$. It is easy to see that
this element depends on the image of $c$ in $H^1(G,\BZ)$ and gives rise
to a map
$\BZ[H^1(G,\BZ)] \to \La_G$. It is easy to see that this map is an 
isomorphism.
\end{proof}

Let us define a vector space
$$
\A^{RV}(\La)=
\left(\oplus_{G} \BZ[H_1(G,\BZ)] \cdot G \right)/(\AS,\IHX, \Aut)
$$
where the sum is over trivalent graphs with oriented edges.

Lemma \ref{lem.Aid1} and \ref{lem.RV} imply that

\begin{corollary}
\lbl{cor.3views}
There are canonical isomorphisms:
$$
\A(\La) \cong 
\A^C(\La)
\cong
\A^{RV}(\La).
$$
\end{corollary}

The map 
$$
\A^{RV}(\La) \to \G \SO
$$ 
(which is the composite of \eqref{eq.ARla} with the identification of the 
above corollary), can be defined more explicitly without appealing 
to finger moves. Indeed, observe that given an embedding $\vphi:
G \to S^3\sminus \O$, then taking linking number in $S^3$
of $\O$ with $1$-cycles in $\vphi(G)$, we get a canonical element 
$\vphi^\ast(\O) \in H^1(G,\BZ)$. 

\begin{lemma}
\lbl{lem.lf3}
If $\vphi, \psi:G\to S^3\sminus \O$ are two
embeddings of a trivalent graph of degree
$2n$ such that $\vphi^\ast(\O)=\psi^{\ast}(\O) \in H^1(G,\BZ)$,
then $[\SO,\vphi(G)]=[\SO,\psi(G)] \in \Gl {2n}\SO$.
\end{lemma}

\begin{proof}
Our Sliding Lemma implies that an edge of $G$ can slide past two arcs
of $K$ with opposite orientations. It can also slide past another edge.
It is easy to see that this implies our result.
\end{proof}

Using the above lemma, and given a vertex-oriented trivalent graph
$G$ and an element $c \in H^1(G,\BZ)$, let $\vphi: G \to S^3\sminus\O$ be any
embedding with associated element $c$. This defines a map
$\BZ[H^1(G,\BZ)] \cdot G \to \G \SO$. The $\AS, \IHX$ and $\Aut$ relations
are satisfied and this gives rise to a well-defined map
$$
\A^{RV}(\La) \longto \G \SO.
$$

\begin{remark}
\lbl{rem.vog}
Given a cocommutative Hopf algebra $H$ and an arbitrary graph $G$ with 
oriented edges, Vogel associates an abelian group $\La^V_G$. 
In case $H=\BZ[t,t^{-1}]$ with comultiplication $\Delta(t)=t \otimes t$ and 
antipode $s(t)=t^{-1}$ the abelian group $\La^V_G$ is the one given by 
Definition \ref{def.Lahalf}.
\end{remark}

\begin{remark}
\lbl{rem.allwork}
All the results of this section work without change if we replace
the pair $\SO$ with a pair $(M,K)$ of a knot $K$ in a \ihs\ $M$
with trivial Alexander polynomial.
\end{remark}

\section{The behavior of the Aarhus integral under surgery on claspers}
\lbl{sub.Zbrief}

In this section, we specialize the general principle of defining/calculating
the Aarhus invariant to the case of links obtained by surgery
on claspers. The reader is referred to \cite{A1} for a detailed discussion
of the Aarhus integral. 
In the present paper we will only be interested in pairs $(M,K)$
of knots $K$ in \ihs s $M$, and $(N,\emptyset)$ of \ihs s $N$. In case
$M=S^3$, we will be using the term ``Kontsevich integral'' rather than
LMO invariant or Aarhus integral. Hopefully this will not cause any confusion 
or historical misunderstanding. 

\subsection{A brief review of the Aarhus integral}
\lbl{sub.LMObrief}

Given a nondegenerate framed link $L$ in a $S^3$ 
(i.e., whose linking matrix is invertible over $\BQ$) the Aarhus invariant 
$Z(S^3_L,\emptyset)$ of $S^3_L$ is obtained from the Kontsevich integral
$Z(S^3,L)$ in the following way:
\begin{itemize}
\item
Consider $Z(S^3,L)$, an element of the completed (with respect to the 
Vassiliev degree) vector space $\A(\tcircle_L)$ of chord diagrams on 
$L$-colored disjoint circles.
\item
After some {\em suitable basing} of $L$ (defined below)
we can lift $Z(S^3,L)$ to an element of the completed algebra
$\A(\tline_L)$ of unitrivalent graphs on $L$-colored vertical segments.
\item
Symmetrize the legs on each $L$-colored segment to convert
$Z(S^3,L)$ to an element of the completed algebra $\A(\ast_L)$ of 
unitrivalent graphs with symmetric $L$-colored legs.
\item
Separate the {\em strut} part $Z^q(S^3,L)$
from the other part $Z^t(S^3,L)$:
$$
Z(S^3,L)=Z^q(S^3,L) \, Z^t(S^3,L)
$$
where $Z^t$ contains no diagrams that contain an $L$-labeled strut
component. It turns out that the strut-part
is related to the linking matrix of $L$ as follows:
$$
Z^q(S^3,L)=\exp\left(\frac{1}{2} \sum_{x,y \in L} l_{xy} \st{x}{y}  \right)
$$
\item
Glue the $L$-colored legs of the graphs of $Z^t(S^3,L)$ pairwise in all
possible ways and multiply the created edges by the entries of the negative 
inverse linking matrix of $L$.
\item
Finally renormalize by a factor that depends on the signature of the linking
matrix of $L$. The end result, $Z(S^3_L,\emptyset)$ depends only on
the 3-manifold $S^3_L$ and not on the surgery presentation $L$ of it
or the basing of $L$.
\end{itemize}

What is a ``suitable basing of $L$?'' Ideally, we wish we could choose a base
point on each component of $L$ to convert $L$ into a union of $L$-labeled
intervals. Unfortunately, it is more complicated: a suitable basing of $L$
is a string-link representative of $L$, equipped with {\em relative
scaling} (i.e., a parenthetization) between the strings. Such objects were 
called {\em q-tangles} (see \cite{LM}) by Le-Murakami, {\em non-associative 
tangles} by Bar-Natan (see \cite{B3}), and were described in the equivalent 
language of {\em dotted Morse links} in \cite[part II, Section 3]{A1}.

Notice that certain parts (such as the $L$-labeled struts) of the Kontsevich 
integral of a based link $L$ are independent of the basing.

\begin{remark}
\lbl{rem.Aqhs}
The above discussion of the Aarhus integral also works when we start from
pairs $(M,L)$ of nondegenerate framed links $L$ in \qhs\ $M$.

The discussion also works in a relative case when we start from pairs
$(M, L \cup L')$ (with $L$ nondegenerate in a \qhs\ $M$) and do surgery on  
$L$ alone. 
\qed
\end{remark}

In what follows, we will assume silently that the links in question
are suitably based, using parenthesized string-link representatives.

\subsection{The Aarhus integral and surgery on claspers}
\lbl{sub.LMOy}

We will be interested in links that describe surgery on claspers. 
Consider a clasper $G$ of degree $1$ in $S^3$ (a so-called \ygraph ). Surgery 
on the clasper
$G$ can be described by surgery on a six component link $E \cup L$ associated 
to $G$, where $E$ (resp. $L$) is the three component link that consists of 
the edges (resp. leaves) of $G$. The linking matrix of $E\cup L$
and its negative inverse are given as follows:
\begin{equation}
\lbl{eq.arms}
\left(\begin{array}{cc}
0 & I \\ I & \lk(L_i,L_j) \\
\end{array}\right) \hspace{0.3cm}
\text{ and } \hspace{0.3cm}
\left(\begin{array}{cc}
\lk(L_i,L_j) & -I \\ -I & 0 \\
\end{array}\right) .
\end{equation}
The six component link $E \cup L$ is partitioned in three blocks
of two component links $A_i=\{E_i,L_i\}$ each for $i=1,2,3$, the {\em arms}
of $G$. A key feature of surgery on $G$ is the fact that surgery on any 
proper subset of the set of arms does not alter $M$. In other words,
alternating $M$ with respect to surgery on $G$ equals to alternating $M$ with 
respect to surgery on all nine subsets of the set of arms 
$A=\{A_1,A_2,A_3 \}$. That is,

\begin{equation}
\lbl{eq.altZ}
Z([S^3,G])=Z([S^3, A]).
\end{equation}

Due to the locality property of the Kontsevich integral (explained in 
\cite{L} and in a leisure way in \cite[II, Section 4.2]{A1}), the nontrivial 
contributions to the right hand side of Equation \eqref{eq.altZ} come from
the part of $Z^t(S^3,A_1 \cup A_2 \cup A_3)$ that consists of graphs with
legs {\em touch} (i.e., are colored by) all three arms of $G$.

The above discussion generalizes to the case of an arbitrary 
disjoint union of claspers. 

\subsection{The Aarhus integral is a universal invariant of \ihs s}
\lbl{sub.Auni}

In this section we will explain briefly why the Aarhus integral (evaluated
at \ihs s and graded by the Euler degree) is a universal \fti . We will
use essentially the same ideas to prove Theorem \ref{thm.loop2} in the
following section. 
The ideas are well-known and involve elementary counting arguments, see
\cite{A1,L} and also \cite{Ga2}.

We begin by showing the following well-known proposition, which we could not
find in the literature.

\begin{proposition}
\lbl{prop.MKu1}
The Euler degree $n$ part of $Z(\cdot, \emptyset)$ is a type $n$ invariant 
of \ihs s with values in $\A_n(\phi)$.
\end{proposition}

\begin{proof}
Since $\A_{\text{odd}}(\phi)=0$, it suffices to consider the case of 
even $n$.

Recalling Definition \ref{def.ihsfti}, suppose that $G=\{G_1, \dots, G_{m} \}$ 
(for $m \geq 2n+1$) is a collection of claspers in $S^3$ each of 
degree $1$, and let $A$ denote the set of arms of $G$.
Equation \eqref{eq.altZ} and its following discussion implies that
$$
Z([S^3,G])=Z([S^3, A])
$$
and that the nonzero contribution to the right hand side come from
diagrams in $Z^t(S^3,A)$ that touch all arms. Thus, contributing
diagrams have at least $3(2n+1)=6n+3$ $A$-colored legs, to be glued 
pairwise. Since pairwise gluing needs an even number number of univalent
vertices, it follows that we need at least $6n+4$ $A$-colored legs.

Notice that $Z^t(S^3,A)$ contains no struts. Thus, at most three $A$-colored
legs meet at a vertex, and after gluing the $A$-colored legs we obtain
trivalent graphs with at least $(6n+4)/3=2n+4/3$ trivalent vertices, in other
words of Euler degree at least $2n+2$. Thus,
$Z_{2n}([S^3,G])=0 \in \A_{2n}(\phi)$, which implies that
$Z_{2n}$ is a invariant of \ihs s of type $2n$ with values in $\A_{2n}(\phi)$.
\end{proof}

Sometimes the above vanishing statement is called {\em counting above
the critical degree}. Our next statement can be considered as {\em counting on the critical degree}. We need a preliminary definition. 

\begin{definition}
\lbl{def.cc}
Consider a clasper $G$ in $S^3$ of degree $2n$, and let 
$G^{break}=\{G_1,\dots, G_{2n}\}$
denote the collection of degree $1$ claspers $G_i$ which are obtained
by inserting a Hopf link in the edges of $G$. Then the 
{\em complete contraction} $\la G\ra\in\A(\phi)$ of $G$ is defined to be the 
sum over all ways of gluing pairwise the legs of $G^{break}$, where we multiply
the resulting elements of $\A(\phi)$ by the product of the linking numbers of 
the contracted leaves. 
\end{definition}

\begin{proposition}
\lbl{prop.MKu2}
If $G$ is a clasper of degree $2n$ in $S^3$, then
$$
Z_{2n}([S^3,G^{break}])=\la G \ra \in \A_{2n}(\phi)
$$
\end{proposition}
 
\begin{proof}
It suffices to consider a collection
$G=\{G_1, \dots, G_{2n} \}$ of claspers in $S^3$ each of 
degree $1$. Let $A$ denote the set of arms of $G$. The counting argument
of the above Proposition shows that the contributions
to $Z_{2n}([S^3,G])=Z_{2n}([S^3, A])$ come from complete contractions
of a disjoint union  $D=Y_1 \cup \dots \cup Y_{2n}$ of $2n$ {\em vortices}. 
A vortex is the diagram $\sfY$, the next simplest unitrivalent graph after
the strut. Furthermore, the $6n$ legs of $D$ should touch all
$6n$ arms of $G$. In other words, there is a 1-1 correspondence between
the legs of such $D$ and the arms of $G$.

Consider a leg $l$ of $D$ that touches an arm $A_l=\{E_l,L_l\}$ of $G$.
If $l$ touches $L_l$, then due to the restriction of the negative inverse
linking matrix of $G$ (see Section \ref{sub.LMOy}), it needs to be contracted
to another leg of $D$ that touches $E_l$. But this is impossible, since
the legs of $D$ are in 1-1 correspondence with the arms of $G$.

Thus, each leg of $D$ touches presicely one edge of $G$. In particular, each 
component $Y_i$ of $D$ is colored by three edges of $G$. 

Notice that the set of edges of $G$ is an algebraically split link.
Given a vortex colored by three edges of $G$, the coefficient of it
in the Kontsevich integral equals to the {\em triple Milnor invariant}, 
(as is easy to show, see for example \cite{HM}) and vanishes unless all three 
edges are part of a degree $1$ clasper $G_i$. When the triple Milnor invariant
does not vanish, it equals to $1$ using the orientation of the clasper $G_i$.

Thus, the diagrams $D$ that contribute are a disjoint union of $2n$ vortices 
$Y=\{Y_1,\dots,Y_{2n}\}$ and these vortices are in 1-1 correspondence with
the set of claspers $\{G_1,\dots,G_{2n}\}$, in such a way that the
legs of each vortex $Y_i$ are colored by the edges of a unique clasper $G_j$.

After we glue the legs of such $Y$ using the negative inverse linking matrix 
of $G$, the result follows.
\end{proof}

\begin{remark}
\lbl{rem.picky}
Let us mention that the discussion of Propositions \ref{prop.MKu1} and 
\ref{prop.MKu2} uses the
unnormalized Aarhus integral; however since we are counting above the critical
degree, we need only use the degree $0$ part of the normalization which
equals to $1$; in other words we can forget about the normalization.
\end{remark}

The above proposition is useful in realization properties of the $Z_{2n}$
invariant, but also in proving the following {\em Universal Property}:

\begin{proposition}
\lbl{prop.MKu3}
For all $n$, 
the composite map of Equation \eqref{eq.AS3} and Proposition \ref{prop.MKu2}
$$
\A_{2n}(\phi)\longto \G_{2n}(S^3) \stackrel{Z_{2n}}\longto \A_{2n}(\phi)
$$
is the identity. Since the map on the left is onto, it follows that
the map \eqref{eq.AS3} is an isomorphism, over $\BQ$. 
\end{proposition}

We may call the map $G \to Z_{2n}([S^3,G])$ (where $G$ is of degree $2n$)
the {\em symbol} of $Z_{2n}$.

\begin{remark}
In the above propositions \ref{prop.MKu1}-\ref{prop.MKu3}, $S^3$ can be 
replaced by any \ihs\ (or even a \qhs s) $M$.  \qed
\end{remark}

\subsection{Proof of Theorem \ref{thm.loop2}}
\lbl{sub.thm.loop2}

Our goal in this section is to modify (when needed) the discussion of 
Sections \ref{sub.LMObrief}-\ref{sub.Auni} and deduce in a natural way
the proof of Theorem \ref{thm.loop2}. It is clear that struts and vortices
and their coefficients in the Kontsevich integral play an important role.

As we will see presently, the same holds here, when we replace struts and
vortices by their {\em hairy} analogues:

\begin{figure}[htpb]
$$ 
\psdraw{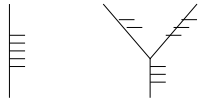}{1in}
$$
\caption{Hairy struts and hairy vortices.}\lbl{HH}
\end{figure}

We begin by considering an $\SO$-null clasper $G$. We know already that
we can compute $Z([S^3,G])$ from the Kontsevich integral $Z(S^3, G \cup \O)$
of the corresponding link by integrating along $G$-colored struts.
Choose string-link representatives of $G \cup \O$ with relative scaling. 
We will denote edges and leaves of $G$ by $E$ and $L$ respectively.

The diagrams that appear in $Z(S^3,G \cup \O)$ are {\em hairy}, in the sense
that upon removal of all $\O$-colored legs, they are simply unitrivalent 
graphs whose legs are colored by $G$. 

Let us separate the hairy strut part $Z^{hq}(S^3,G \cup \O)$ from the other 
part $Z^{ht}(S^3, G \cup \O)$:
$$
Z(S^3,G \cup \O)=Z^{hq}(S^3,G \cup \O) \, Z^{ht}(S^3,G \cup \O)
$$
using the disjoint union multiplication,
where $Z^{ht}$ contains no diagrams that contain an $L$-labeled hairy strut
component. Let us write the hairy strut part as follows: 
$$
Z^{hq}(S^3, G \cup \O)=\exp\left(\frac{1}{2} \sum_{x,y \in L \cup E} 
f(L_x,L_y,\O) \right)
$$
where $E \cup L$ is the corresponding link of leaves and edges of $G$ and 
$$
f(L_x,L_y,\O)=\sum_{n=0}^\infty \mu_{L_x,L_y,O \dots n \, 
\text{times} \dots 0} 
\psdraw{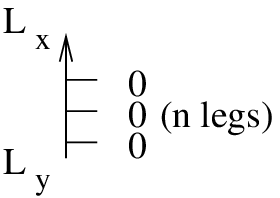}{1in}
$$
In the above equation, $\mu_{L_x,L_y,O \dots n \, 
\text{times} \dots 0}$ is the coefficient in $f(L_x,L_y,\O)$ of the
hairy diagram appearing on the right.
We will call $f( \cdot, \cdot, \O)$ the {\em hairy linking matrix}.

\begin{lemma}
\lbl{lem.hairvanish}
The hairy linking matrix of $G$ and its negative inverse are given by 
$$
\left(\begin{array}{cc}
0 & I \\ I & f(L_i,L_j,\O) \\
\end{array}\right) \hspace{0.3cm}
\text{ and } \hspace{0.3cm}
\left(\begin{array}{cc}
f(L_i,L_j,\O) & -I \\ -I & 0 \\
\end{array}\right) .
$$
\end{lemma}

\begin{proof}
Since $\{E_i,E_j,\O\}$ is an unlink and $\{E_i,L_{j \neq i},\O\}$ 
is the disjoint union of an unknot $E_i$ and the link $\{L_{j \neq i},\O\}$,
it suffices to consider the case of a hairy strut $\sfI^{E_i}_{L_i}$. 
In this case,
$E_i$ is a meridian of $L_i$. The formula for the Kontsevich integral of 
the Long Hopf Link of \cite{BLT} (applied to $E_i$), together with the fact 
that $L_i$ has linking number zero with $\O$, imply that the hairy part 
$\sfI^{E_i}_{L_i}$ vanishes. 
\end{proof}

\begin{lemma}
\lbl{lem.f}
$Z([\SO,G])\in \A(\ast_\O)$ can be computed from $Z^{ht}(S^3,G\cup \O) 
\in \A(\ast_{G \cup \O})$ by gluing pairwise the legs of the $G$-colored 
graphs of $Z^{ht}(S^3,G\cup \O)$ using the negative inverse hairy linking
matrix.
\end{lemma}

\begin{proof}
We have that
$$
Z^t(S^3, G \cup \O)= \exp\left(\frac{1}{2} \sum_{x,y \in L \cup E} 
\st{x}{y} f(x,y,\O)-\lk(x,y)  \right) \, 
Z^{ht}(S^3, G \cup \O) 
$$
In other words, the diagrams that contribute in $Z([\SO,G])$ are those
whose components either lie in $Z^{ht}(S^3, G \cup \O)$ or are hairy
struts of the shape $\sfI^{L_i}_{L_j}$. Because of the restriction 
of the linking matrix \eqref{eq.arms}, the pair of legs $\{L_i,L_j\}$
of a hairy strut of the above shape must be glued to a pair of legs
labeled by $\{E_i,E_j\}$ as follows:
$$
\psdraw{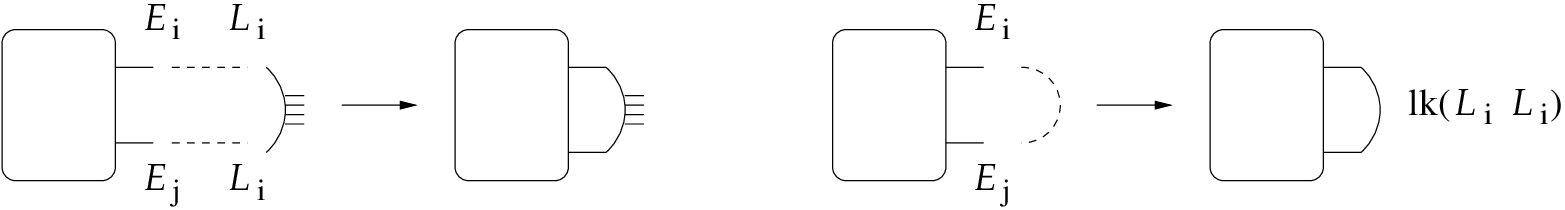}{5in}
$$
The result of this gluing is equivalent to gluing pairs of $(E_i,E_j)$ 
colored legs using the negative inverse hairy linking matrix. 
This concludes the proof of the lemma.
\end{proof}

Given this Lemma, the proof of Proposition \ref{prop.MKu1} works without
change and proves Theorem \ref{thm.loop2}.

Furthermore, observe that the coefficients of hairy vortices with legs colored
by the edges of $G$ and with nonzero number of hair vanish. This is true since
the Borromean rings (i.e., the edges of $G$) form an unlink in the complement
of $\O$. Thus, the proof of Proposition \ref{prop.MKu2} implies the following

\begin{theorem}
\lbl{thm.u2}
If $G$ is a $\SO$-null clasper of degree $2n$, then
$$
Z_{2n}([\SO,G])=\la G \ra \in \A_{2n}(\ast)
$$
where complete contractions are using the matrix of Lemma \ref{lem.f} instead 
and put the entries of it as beads on the edges that are created by the 
contraction of the legs.
\end{theorem}

Let us remark that the analogue of Proposition \ref{prop.MKu3} does not hold
since $Z$ takes values in $\A(\ast)$ rather than in $\A(\La)$. The invariant
$\Zrat$ discussed in Section \ref{sub.further} takes values in $\A(\La)$,
satisfies the analogue of \ref{prop.MKu2} and is thus a universal $\BQ$-valued
invariant of Alexander polynomial $1$ knots, with respect to the null-move.

\section{Abelian invariants: equivariant linking numbers}
\lbl{sec.abelian}

In Theorem \ref{thm.u2} we calculated the symbol $Z_{2n}[\SO,G]
\in \A_{2n}(\La)$ in terms of a complete contraction of $G$ that uses the
the hairy linking matirx. Struts correspond to linking numbers and 
we will show that hairy struts correspond to equivariant linking numbers.
The goal of this section is to discuss the latter.

\subsection{A review of linking numbers}
\lbl{sub.linkingn}

Recall that all links are oriented. We begin by recalling the definition of 
the {\em linking number} $\lk(L)$ of a link $L=(L_1,L_2)$ of two ordered
components in $S^3$. Since $H_1(S^3,\BZ)=0$, there is an
oriented surface $\S_1$ that $L_1$ bounds. We then define 
$\lk(L)=[\S_1]\cdot [L_2] \in \BZ$ where $\cdot$ is the intersection
pairing. Since $H_2(S^3)=0$, the result is independent of the choice of
surface $\S_1$.

The following lemma summarizes the well-known properties of the linking 
numbers of two component links in $S^3$:

\begin{lemma}
\lbl{lem.mu}
{\bf (Symmetry)}
$\lk(L)$ does not depend on the ordering of the components of $L$. \newline
{\bf (Cutting)}
If a component $L_i$ of $L$ is a connected sum of
$L_i'$ and $L_i''$ (as in Figure \ref{split3}),  then
\begin{eqnarray*}
\lk(L)& = &
\lk(L')
+\lk(L'') 
\end{eqnarray*}
{\bf (Initial Condition)}
$$
\lk(\psdraw{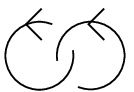}{0.3in})=1, \,\,\,
\lk(\psdraw{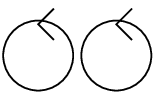}{0.3in})=0. 
$$
{\bf (Uniqueness)}
If a function of two-component links satisfies the Symmetry, Cutting, and 
Initial Conditions, then it equals to $\lk$.
\end{lemma}

\begin{proof}
It is easy to see that linking numbers satisfy the above axioms.

For Uniqueness, assume that $\aa$ is another such function, and
let $\bb=\aa-\lk$. consider a link $L=(L_1,L_2)$, and a surface $\S_1$
that bounds $L_1$. Using the Cutting Property, as in the proof of Theorem
\ref{thm.GGPAS3}, it follows that $\bb(L)$ is a linear combination of 
$\bb(\text{Hopf Link})$ and $\bb(\text{Unlink})$, and hence zero.
\end{proof}

\begin{remark}
\lbl{rem.qlk}
The above definition of linking number can be extended to the case of an 
two component links $L$ in a \qhs\ $M$. In that case, the 
linking number takes values in $1/|H_1(M,\BZ)| \BZ$, and Proposition
\ref{lem.mu} continues to hold. \qed
\end{remark}

\begin{remark}
\lbl{rem.qlk2}
The above definition of linking number can be extended to the case of a 
two-component link $L$ in a possibly open 3-manifold $N$ that satisfies
$H_1(N,\BZ)=H_2(N,\BZ)=0$. \qed
\end{remark}

\subsection{Equivariant linking numbers}
\lbl{sub.linkingf}

In this section we review well-known results about equivariant linking 
numbers. These results are useful to the {\em surgery view} of knots,
(see \cite{L3,KY,Ro} and also \cite{GK1}) and will help us identify the hairy 
struts with equivariant linking numbers.

We begin with the following simple situation: consider a null homologous
link $L=(L_1,L_2)$ of two ordered components in $X=S^3\sminus\O$. 

\begin{definition}
\lbl{def.Lnull}
We will call such links $L$ $(S^3,\O)$-{\em null}.
\end{definition}

Consider the
universal abelian cover $\ti X\to X$ corresponding to the natural map 
$\pi_1(X) \to H_1(X,\BZ)\cong \BZ$. Since $L$ is null homologous,
it lifts to a link $\ti L$, which is invariant under the action by
the group of deck transformations $\BZ$.  Note that $\ti X$ is an open
3-manifold diffeomorphic to $D^2 \times \BR$, in particular $H_1(\ti X,
\BZ)=H_2(\ti X, \BZ)=0$. Using Remark \ref{rem.qlk2}, we can define
linking numbers between the components of $\ti L$.

Consider an arc-basing $\ga$ of $L$, that is an embedded arc in $X$ that 
begins in one component of $L$ and ends at the other, and is otherwise 
disjoint from $L$. Then, we can consider a lift of $L \cup \ga$ to $\ti X$
which is an arc-based two component link $(\ti L_1, \ti L_2)$ in $\ti X$.

\begin{definition}
\lbl{def.lkf}
If $(L,\ga)$ is an $\SO$-null arc-based link, we define the {\em
equivariant linking number} by
$$
\lgt(L)=\sum_{n \in \BZ} 
 \lk(\ti L_1 , t^{n} \ti L_2 ) \, t^n \in \La
$$
\end{definition}

The above sum is finite. Moreover, since linking numbers (in $\ti X$)
are invariant under the action of deck transformations, it follows that
the above sum is independent of the choice of lift of $L \cup \ga$ to $\ti X$.

The following lemma summarizes the properties of the equivariant linking 
number of $\SO$-null arc-based links. $\e: \La \to \BZ$ stands for the map 
$t \to 1$.

\begin{lemma}
\lbl{lem.lkprops}
{\bf (Symmetry)} If $\s$ is a permutation of the two components of $L$ then,
we have 
$$\lgt(L_{\s})(t)= \lgt(L)(t^{-1})
$$ 
{\bf (Specialization)}
$$\e \,\, \lgt(L)=\lk(L)
$$
{\bf (Cutting)}
Suppose that a component $L_i$ of an $(S^3,\O)$-null link
$(L,\ga)$ is a connected sum of
$L_i'$ and $L_i''$ (as in Figure \ref{split3}), and $(L',\ga)$ and 
$(L'',\ga)$ are $(S^3,\O)$-null. Then,
\begin{eqnarray*}
\lgt(L)& = &
\lgt(L')
+\lgt(L'') 
\end{eqnarray*}
{\bf ($\La$-Sliding)} 
If $(L^s,\ga)$ denote the result of
sliding the arc-basing of the first component of $(L,\ga)$ 
along an oriented arc-based curve $S$, 
then 
\begin{eqnarray*}
\lgt(L^s)& = & 
t^l\lgt(L) 
\end{eqnarray*}
where $l=\lk(S,\O)$. \newline
{\bf (Initial Condition)}
If $L$ lies in a ball disjoint from $\O$, then 
\begin{eqnarray*} 
\lgt(L) &=& \lk(L).
\end{eqnarray*}
{\bf (Uniqueness)} If a function of $\SO$-null arc-based links
satisfies the Symmetry, $\La$-Sliding, Cutting and Initial Conditions,
then it equals to the equivariant linking number.
\end{lemma}

\begin{proof}
The Symmetry, Cutting and $\La$-Sliding Properties follows immediately from 
Lemma \ref{lem.mu}.

The Specialization Property follows from the
well-known fact that given a covering space $\pi: \ti X \to X$ and a cycle
$c$ in $X$ that lifts to $\ti c$ in $\ti X$, and a cycle $c'$ in $\ti X$,
then the intersection of $\ti c$ with $c'$ equals to the intersection of $c$
with the push-forward of $c'$ in $X$.   

The Uniqueness statement follows from the proof of Theorem \ref{thm.loop1}
given in Section \ref{sub.graded}.
\end{proof}

\begin{remark}
\lbl{rem.Kojima}
The definition \ref{def.lkf} and Lemma \ref{lem.lkprops} can be extended 
without change to the case of a null homologous arc-based two component 
link $(L,\ga)$ in $M\sminus K$ where $(M,K)$ is a knot $K$ in a \ihs\ $M$ 
with trivial Alexander polynomial. Notice that
$H_2(\ti {M\sminus K},\BZ)=0$ and the condition on the Alexander polynomial
ensures that $H_1(\ti {M\sminus K},\BZ)=0$. 

A further generalization of definition \ref{def.lkf} and Lemma 
\ref{lem.lkprops} is possible to the case of a null homologous arc-based two 
component link $(L,\ga)$ in $M\sminus K$ where $(M,K)$ is a knot $K$ in a 
\ihs\ $M$. In that case, $\lgt(L)$ lies in a {\em localization} $\Lloc$ 
of $\La$ defined by
$$
\Lloc=\{ p/q \, | p, q \in \La, \, q(1)=\pm 1 \}
$$
This specializes to tha case of $L$ being a zero-framed parallel of a 
null homologous knot $K'$ in $S^3 \sminus K$ and $\ga$ a short arc between
the two components of $L$. In this case, $\lgt(K',K')$ coincides with
the self-linking $\eta$-function $\eta(K,K')$ 
of Kojima-Nakanishi, \cite{KY}.
\qed
\end{remark}

\section{Hairy struts are equivariant linking numbers}
\lbl{sub.Kbehavior}

The purpose of this section is to show that the hairy struts coincide
with equivariant linking numbers. 

We begin with a definition. A {\em special link} $L \subset S^3$ is 
the union of $\O$ and an $\SO$-null link. Note that a special link
$L$ has a special component, namely $\O$. 

In what follows, $\ga$ will refer to a {\em disk-basing} (i.e., a string-link
representative) of $L$. Consider a 
disk-based special link $(L,\ga)$, together with a choice of relative
scale, and its Kontsevich integral. There is an algebra isomorphism 
$\s: \A(\tline_L)\to\A(\ast_L)$ with inverse the symmetrization
map $\chi: \A(\ast_L) \to \A(\tline_L)$ which is the average of all
ways of placing symmetric $L$-colored legs on $L$-intervals. In what follows, 
we will denote by $Z(S^3,L)$ the image of the Kontsevich integral in either
algebra, hopefully without causing confusion.

Consider a homotopy quotient $\A^h(\ast_{L})$ of the algebra 
$\A(\ast_{L})$ where we quotient out by all nonforests (i.e., graphs
with at least one connected component which is not a forest)
and  by all {\em flavored} 
forests that contain a tree with at least two legs flavored by the same 
component of $L$, see also \cite{B2} and \cite{HM}. A {\em flavoring}
of a graph is a decoration of its univalent vertices by a set of colors,
which in our case is $L$. We will be interested 
in the Lie subalgebra $\A^{c,h}(\ast_{L})$ of $\A^h(\ast_{L})$ spanned by 
all connected diagrams. It is spanned by hairy trees colored by $L \sminus\O$ 
such that each label of $L\sminus \O$ appears at most once. 
Notice that if $T$ and $T'$ are $A$-colored and $A'$-colored hairy trees
then 
$$
[T,T']=
\begin{cases}
0 & \text{if $|A \cup A'|\neq 1$}, \\
T\cdot_{i} T' & \text{if $A \cup A'=\{i\}$}
\end{cases}
$$
where $T\cdot_i T'$ is the result of {\em grafting} the trees $T$ and $T'$
along their common leg $i$. 

Let us define $Z^h(S^3,L)$ to be the image of
$Z(S^3,L)$ in $\A^{h}(\ast_{L})$.
The group-like property of the Kontsevich integral implies that
$\log Z^h(S^3,L)$ lies in $\A^{c,h}(\ast_{L })$.
In the case $L$ has three components, 
$\A^{c,h}(\ast_{L})$ has a basis consisting of the
trees $T_n(12)$, for $n \geq 0$, together with the degree $1$ struts $t_{10}$
and $t_{20}$, and that
\begin{equation}
\lbl{eq.commute}
\ad_{t_{20}}^n t_{12}=(-1)^nT_n(12) \text{ and }
\ad_{t_{10}}^n t_{12}= T_n(12).
\end{equation}
Here the product (and the commutator) are taken with respect to the
natural multiplication on $\A(\tline_L)$.

\begin{figure}[htpb]
$$ 
\psdraw{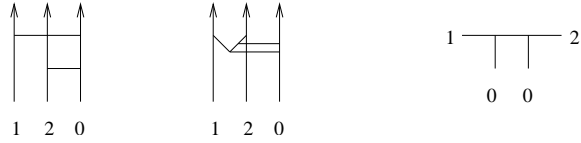}{3.0in} 
$$
\caption{On the left, $t_{20} t_{10}$, where we multiply from left to right 
and from the bottom to the top. In the middle, $T_2(12)$. In the right, an 
alternative view of $T_2(12)$ with symmetrized legs.}
\end{figure}

It follows by the definition of the hairy linking matrix that for every 
special link $L$, we have
\begin{equation}
\lbl{eq.Zh}
\log Z^h(S^3,L)=\sum_{L'}
\fg(L')
 \in \A^{c,h}(\ast_{L})
\end{equation}
where we explicitly denote the dependence on the disk basing $\ga$, and where
the sum is over all special sublinks $L'$ of $L$ that contain $\O$
and two more components of $L\sminus \O$ (with repetition).

At this point, let us convert hairy struts into a power series as follows:

\begin{definition}
\lbl{def.pseries}
If $L$ is a special link of three components with disk-basing $\ga$ and
choice of relative scaling, and
$$
f(L)=\sum_{n=0}^\infty \mu_{L_x,L_y,O \dots n \, 
\text{times} \dots 0} 
\psdraw{attach2}{1in}
$$ 
then
$$
\phig(L)=\sum_{n=0}^\infty \mu_{L_x,L_y,O \dots n \, 
\text{times} \dots 0} \, x^n \in \BQ[[x]]
$$
\end{definition}

\begin{lemma}
\lbl{lem.paren}
If $(L,\ga)$ is a special link of three components equipped with
relative scaling, then $\phig(L)$ is independent of the relative scaling of 
$L$. 
\end{lemma}

\begin{proof}
It will be more convenient to present the proof in the algebra
$\A(\tline_L)$.
Let $\s'$ denote a change of scaling of a string-link $\s$, as shown
in the following figure in case $\s$ has three stands:
$$\psdraw{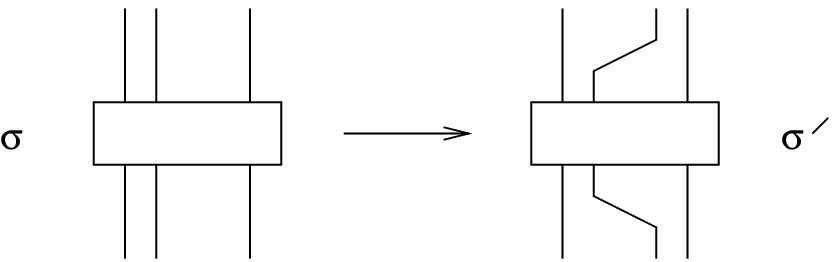}{2in}$$
The locality of the Kontsevich integral implies that $Z(\s')=\Phi Z(\s)
\Phi^{-1}$ for an associator $\Phi$. Write $\Phi=e^{\phi}$ for 
an element $\phi \in [\cL,\cL]$, where
$\cL$ is the free-Lie algebra $\cL$ of two generators $a,b$.
The following identity 
\begin{equation}
\lbl{eq.flie}
e^{a} e^b e^{-a}=e^{\exp(\ad_a) b}
\end{equation}
(valid in a free Lie algebra of two generators)
implies that $\log Z(\s')= e^{\ad_{\phi}} \log Z(\s) \in \A(\tline_{123})$.
If we project the above equality to the quotient $\A^h(\tline_{123})$
where connected diagrams with two legs either both on the first strand
or both on the second strand vanish, then it follows that
we can replace $\phi \in [ \cL, \cL ]$ by its image $\bar\phi \in [ \cL, 
\cL]/[\cL, [\cL,\cL]]$. It is easy to see that $\bar\phi=1/24 \, [a,b]$
for {\em any} associator $\Phi$.

Now we can finish the proof of the lemma as follows.
Consider a sting-link $\s$ with relative scaling
obtained from a disk-basing of a special link
$L$ of three components and let $L'$ be the one obtained by a change of 
relative scaling of $L$. Projecting to $\A^h(\tline_{L})$, and  using the 
fact that $\log Z^h(S^3,L)$ lies in the center of the Lie algebra
$\A^{c,h}(\tline_{L})$, it follows that $\log Z^h(S^3,L' )=\log Z^h(S^3,L)$. 
\end{proof}


\begin{lemma}
\lbl{lem.fcut}
$\phig$ satisfies the Cutting Property of Lemma \ref{lem.lkprops}
with $t=e^{x}$.
\end{lemma}

\begin{proof}
Consider the special link $(L_{1'1''20},\ga)$, (this is an abbreviation
for 
$((L_{1'},L_{1''},L_2,\O),\ga)$) whose connected sum
of the first two components gives $(L,\ga)$.
How does the Kontsevich integral of $L_{1'1''20}$ determine
that of $L_{120}$? The answer, though a bit complicated,  is known
by \cite[Part II, Proposition 5.4]{A1}. Following that notation,  we have
$$ Z(S^3,L_{120})=\la \exp(\La^{1'1''}_{1}), 
Z(S^3,L_{1'1''20}) \ra_{1',1''} \in
\A(\ast_{L_{120}}),
$$
where $\la A,B \ra_{1',1''}$ is the operation that glues all $\{1',
1''\}$-colored
legs of $A$ to those of $B$ (assuming that the number of legs of color $1'$
and of color $1''$ in $A$ and $B$ match; otherwise it is defined to be zero),
and 
$$  \La^{1'1''}_{1}= |^{1'}_{1} + |^{1''}_{1} + 
\La^{1'1''}_{1}\text{(other)}
$$
where $\La^{1'1''}_{1}\text{(other)}$
is an (infinite) linear combination of 
of rooted trees with at least one trivalent vertex whose
leaves are colored by $(L_1',L_1'')$ and whose root is colored by $L_1$. 
We will call such trees $(1',1''; 1)$-trees. The reader may consult
\cite[Part II, Proposition 5.4]{A1} for the first few terms of
$\La^{1'1''}_{1}$, which are given by {\em any} Baker-Cambell-Hausdorff 
formula, translated in terms of rooted trees. 

Upon projecting the answer to the quotient $\A^{h}(\ast_{L_{120}})$,
the above formula simplifies. Indeed, if we glue some disjoint union of
trees of type $(1',1'';1)$ that contain at least one trivalent vertex
to some hairy $L_{1'1''2}$-colored trees,
the resulting connected graph will either have nontrivial homology, or at least
two labels of $L_2$. Such graphs vanish in $\A^{h}(\ast_{L_{120}})$.
Thus, when projecting the above formula to $\A^{h}(\ast_{L_{120}})$,
we can assume that $\La^{1'1''}_{1}= |^{1'}_{1} + |^{1''}_{1}$.
Using the fact of how the Kontsevich integral of a link determines that
of its sublinks and the above, it follows that
\begin{eqnarray*}
\log Z^h(S^3,L_{120}) & = & \log \la \exp(|^{1'}_{1} + |^{1''}_{1}),
Z^h(S^3,L_{1'1''20}) \ra_{1',1''} \\
& = &
\log Z^h(S^3,L_{1'20}) + \log Z^h(L_{1''20})
\end{eqnarray*}
which, together with Equation \eqref{eq.Zh} concludes the proof.
\end{proof}

\begin{lemma}
\lbl{lem.fslide}
$\phig$ satisfies the $\La$-Sliding
Property of Lemma \ref{lem.lkprops} with $t=e^{x}$.
\end{lemma}

\begin{proof}
Consider the link $(L_{S120},\ga)$. Recall that a slide move is given by:
\begin{figure}[htpb]
$$
\psdraw{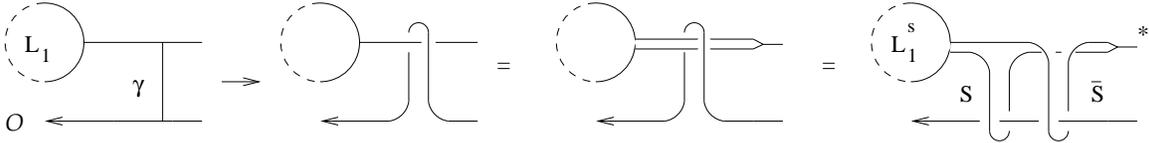}{6in}
$$
\caption{A finger move along an arc $\ga$.}\lbl{finger}
\end{figure}

In an artistic way, the next figure shows the result of a slide move
(compare also with Figure \ref{finger}) that replaces 
$L_1$ 
by $L_1^s:=S \sharp L_1 \sharp \overline{S}$, where $\overline{S}$
is the orientation reversed knot.
$$\psdraw{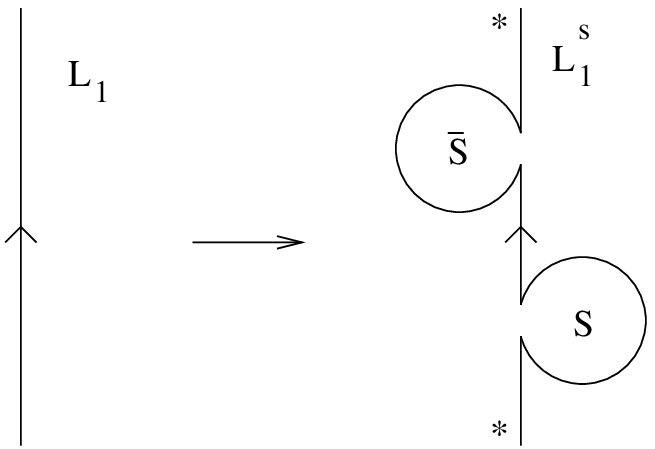}{1.7in}
$$
As in the previous lemma, we have that
$$
Z(S^3,L_{1^s20})=\la \exp(\La^{S1}_{1^s}), \D_{S\oS}Z(S^3,L_{S120})
\exp(\La^{1\oS}_{1^s}) \ra_{1,S,\oS} \in \A(\ast_{L_{1^s20}})
$$
where $\D_{S\oS}$ is the operation that replaces an $S$-colored leg to
an $S+\oS$-colored one.

Upon projecting to $\A^h(\ast_{L_{1^s20}})$, the above formula simplifies.
Indeed, the $\{S,L_1,L_2\}$-colored hairy trees are (after removal of
the hair) of two shapes
$\sfY$ and $\sfI$. 
When we glue, a $\{S,L_1,L_2\}$-colored hairy
$\mathsf{Y}$ vanishes in $\A^h(\ast_{L_{1^s20}})$.
The remaining $\{S,L_1,L_2\}$-colored hairy trees are of the shapes
$
\sfI^{S}_{L_1}, \sfI^{S}_{L_2}, \sfI^{L_1}_{L_2}, 
\sfI^{\oS}_{L_1}, \sfI^{\oS}_{L_2},
$
as well as the hairless $\sfI^{S}_{\O}$ and $\sfI^{\oS}_{\O}$.
When glued to $(S,1; 1^s)$ and $(\oS,1; 1^s)$ rooted trees,
the ones of shapes $
\sfI^{S}_{L_1},  
\sfI^{\oS}_{L_1}$ vanish in  $\A^h(\ast_{L_1^s,L_2,\O})$, thus
we remain it remains to consider only trees of shape 
$\sfI^{S}_{L_2}, \sfI^{L_1}_{L_2}, \sfI^{\oS}_{L_2}$ and the hairless
$\sfI^{S}_{\O}$ and $\sfI^{\oS}_{\O}$.
When we glue these to $(S,1; 1^s)$ and $(\oS,1; 1^s)$ rooted trees,
the only nonzero contribution comes from gluings like 
$$\psdraw{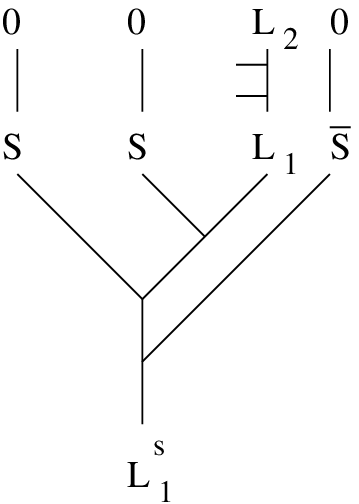}{1.0in}=l^2(-l) \psdraw{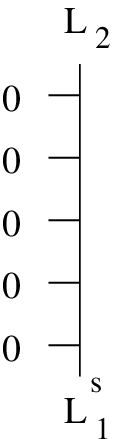}{0.4in}$$
The above figure shows that each of the above gluings can be thought of
as starting from a hairy graph of shape
$\sfI^{L_1}_{L_2}$, and adding to it some additional hair, first on the left
and then on the right; each time
multiplying the result by $l^n(-l)^m$ (where $l=\lk_M(S,\O)$) and $n,m$ is the
number of left and right added hair.
On the other hand, adding hair is the same as commuting with $t_{10}$
(as follows by Equation \eqref{eq.commute}). Translating from gluings back
to Lie algebras, it follows that
\begin{eqnarray*}
Z^h(S^3,L_{1^s20}) & = & \left(e^{\sfI^K_S}
Z^h(S^3,L_{120}) e^{\sfI^K_{\oS}}\right)/(S,\oS \to 1)/(1\to 1^s) \\
& = & 
\left(e^{lt_{10}}
Z^h(S^3,L_{120}) e^{-lt_{10}}\right)/(1\to 1^s)
\end{eqnarray*}
where $E/(1\to 1^s)$ means to replace the label $1$ by $1^s$ in the expression
$E$.
Equation \eqref{eq.flie} implies that
$$
\log Z^h(S^3,L_{1^s20}) =
e^{l \ad_{t_{10}}}
\fg(L_{12})(\ad_{t_{10}})t_{12}/(1\to 1^s) 
$$
which concludes the proof.
\end{proof}

The following theorem is the main result of this section.

\begin{theorem}
\lbl{thm.flk}
For an $\SO$-null arc-based link $(L,\ga)$ of three components
we have that
$$\phig(L)(x)=\lgt(L\sminus\O)(e^{x}) \in \BQ[e^{\pm x}].
$$
\end{theorem}

In other words, the coefficients of hairy struts are equivariant linking
numbers, in particular they are Laurent polynomials.
 
\begin{proof}
This follows from Lemma \ref{lem.lkprops} once we show that
$\phig$ satisfies the Symmetry, Specialization, $\La$-Sliding, 
Cutting and Initial Condition stated in that lemma.

The symmetry follows by Equation \eqref{eq.commute}.
Specialization follows from the fact that $\phig(L)(0)$ coincides
with the coefficient of a strut in the Kontsevich integral, which equals to
$\lk(L_1,L_2)$. The $\La$-Sliding property follows from Lemma 
\ref{lem.fslide}, the Cutting property follows from Lemma \ref{lem.fcut}.
The Initial Condition property follows from the fact that the Kontsevich
integral is multiplicative for disjoint union of links, thus the only
diagrams $T_n$ that contribute to the Kontsevich integral (and also
in $\phig$) in this case is $T_0$, which contributes the linking
number of $L_1$ and $L_2$.
\end{proof}

\begin{remark}
\lbl{rem.MKf}
Theorem \ref{thm.flk} remains true for $(M,K)$-null two component links
$L$ where $K$ is a knot in an \ihs . \qed
\end{remark}

Since the coefficient $\mu_n$ of $x^n$ in the power series
$\phig(L)$ are given by a combination of Milnor's invariants
(as follows by the work of Habegger-Masbaum \cite{HM} for $M=S^3$
and forthcoming work of the first author for a general \ihs ), 
we obtain that

\begin{corollary}
\lbl{cor.conc}
$\lgt$ is a concordance invariant of
$L$ and a link homotopy invariant of $L\sminus K$.
\end{corollary}

This generalizes a result of Cochran, \cite{Co}, who showed
that a power series expansion of the Kojima self-linking function 
was a generating function for a special class of 
repeated Milnor invariants of type $\mu_{11KKKK}$ of a two component 
algebraically split link $(L_1,K)$. 

\subsection{Addendum}
\lbl{sub.addendum}

In an earlier version of the paper, we also identified hairy struts with
equivariant triple Milnor linking numbers, for special links of four 
components. Using this, we gave an algorithm for computing the $2$-loop
part (or, equivalently, the Euler degree $2$ part) of the Kontsevich integral 
of a knot with trivial Alexander polynomial.
This involves a more delicate counting below the critical degree, in the
language of Section \ref{sub.LMOy}. For an easier digestion of the present
paper, we prefer to come back to this matter in a future publication.

{\bf Acknowledgement} We wish to thank D. Bar-Natan, J. Conant,  TTQ Le,
J. Levine, P. Teichner and D. Thurston for many stimulating conversations
and the anonymous referee for pointing out several improvements of an earlier
version.

\ifx\undefined\bysame
	\newcommand{\bysame}{\leavevmode\hbox to3em{\hrulefill}\,}
\fi

\end{document}